\documentclass[12pt]{article}
  \usepackage{amsfonts}
  \usepackage[square,sort&compress,comma,numbers]{natbib} 
   \def\R{\mathbb{R}}

   \def\1{{\rm I\mskip -10.5mu 1}} 
   
   \def\e{{\varepsilon}}

   \def\cC{{\cal C}}

   \def\dist{\mathop{\rm dist}\nolimits}
   
   \def\div{\mathop{\rm div}\nolimits}

   \def\no{\noindent}
   \def\proof{\mbox {{\underline {\sf Proof}} \hspace{2mm}}}
   \def\qed{{\hfill {\em q.e.d.}\\\vspace{1mm}}}
   
   \newcommand{\beq}{\begin{equation}}
   \newcommand{\eeq}{\end{equation}}

\newtheorem{df}{Definition}[section]
\newtheorem{prop}[df]{Proposition}
\newtheorem{lemma}[df]{Lemma}
\newtheorem{teo}[df]{Theorem}
\newtheorem{rem}[df]{Remark}

\newtheorem{ex}[df]{Example}

\newtheorem{cor}[df]{Corollary}

 \newcommand{\sezione}[1]{\section{#1}\setcounter{equation}{0}}

\begin{document}


   \title{Nonexistence of solutions for Dirichlet problems with
  supercritical growth in tubular domains 
}


  \maketitle


 \vspace{5mm}

\begin{center}

{ {\bf Riccardo MOLLE$^a$,\quad Donato PASSASEO$^b$}}

\vspace{5mm}

{\em
${\phantom{1}}^a$Dipartimento di Matematica,
Universit\`a di Roma ``Tor Vergata'',\linebreak
Via della Ricerca Scientifica n. 1,
00133 Roma, Italy.\\
e-mail: molle@mat.uniroma2.it}

\vspace{2mm}

{\em
${\phantom{1}}^b$Dipartimento di Matematica ``E. De Giorgi'',
  Universit\`a di Lecce,\linebreak 
P.O. Box 193, 73100 Lecce, Italy.
}
\end{center}

\vspace{5mm}


{\small {\sc \noindent \ \ Abstract.} - 
We deal with Dirichlet problems of the form
$$
\Delta u+f(u)=0 \mbox{ in }\Omega,\qquad u=0\ \mbox{ on }\partial \Omega
$$
where $\Omega$ is a bounded domain of $\R^n$, $n\ge 3$, and $f$ has
supercritical growth from the viewpoint of Sobolev embedding. 
In particular, we consider the case where $\Omega$ is a tubular domain
$T_\e(\Gamma_k)$ with thickness $\e>0$ and centre $\Gamma_k$, a
$k$-dimensional, smooth, compact submanifold of $\R^n$. 
Our main result concerns the case where $k=1$ and $\Gamma_k$ is
contractible in itself. 
In this case we prove that the problem does not have nontrivial
solutions for $\e>0$ small enough. 
When $k\ge 2$ or $\Gamma_k$ is noncontractible in itself we obtain
weaker nonexistence results. 
Some examples show that all these results are sharp for what concerns
the assumptions on $k$ and $f$.  

\vspace{3mm}



\vspace{1mm}

}


\sezione{Introduction}


The results we present in this paper are concerned with existence or
nonexistence of nontrivial solutions for Dirichlet problems of the
form
\beq
\label{*}
\Delta u+f(u)=0\ \mbox{ in }\ \Omega,\qquad u=0\ \mbox{ on
}\ \partial\Omega,
\eeq
where $\Omega$ is a bounded domain of $\R^n$, $n\ge 3$ and $f$ has
supercritical growth from the viewpoint of the Sobolev embedding.

\no Let us consider, for example, the case where $f(t)=|t|^{p-2}t$ $\forall
t\in\R$ (this function obviously satisfies the condition (\ref{f}) we
use in this paper).
In this case, a well known nonexistence result of Pohozaev (see
\cite{Po}) says that the Dirichlet problem
\beq
\label{p}
\Delta u+|u|^{p-2}u=0\ \mbox{ in }\ \Omega,\qquad u=0\ \mbox{ on
}\ \partial\Omega
\eeq
has only the trivial solution $u\equiv 0$ when $\Omega$ is starshaped
and $p\ge{2n\over n-2}$ (the critical Sobolev exponent).

\no On the other hand, if $\Omega$ is an annulus it is easy to find
infinitely many radial solutions for all $p>1$ (as pointed out by
Kazdan and Werner in  \cite{KW}). 
Thus, it is natural to ask whether or not the nonexistence result of
Pohozaev can be extended to non starshaped domains and the existence
result in the annulus can be extended, for example, to all
noncontractible domains of $\R^n$. 

\no Following some stimulated questions pointed out by Brezis,
Nirenberg, Rabinowitz, etc. (see \cite{B,BN}) many results have been
obtained, relating nonexistence, existence and multiplicity of
nontrivial solutions to the shape of $\Omega$ (see \cite{Di,D88,Pmm89,
P93,P94,Pl92,P4092,Ptmna96,Pd98,MPcvpde06,MPaihp06,MPcras02,MPcras2002},
etc.).

\no In the present paper our aim is to show that, even if the Pohozaev
nonexistence result cannot be extended to all the contractible domains
of $\R^n$, one can prove that there exist contractible non starshaped
domains $\Omega$, which may be very different from the starshaped ones and even
arbitrarily close to noncontractible domains, such that the Dirichlet
problem (\ref{p}) has only the trivial solution $u\equiv 0$ for all
$p>{2n\over n-2}$.

\no In order to construct such domains, we use suitable Pohozaev type
integral identities in tubular domains $\Omega=T_\e(\Gamma_k)$ with
thickness $\e>0$ and centre $\Gamma_k$, where $\Gamma_k$ is a
$k$-dimensional, compact, smooth submanifold of $\R^n$.

\no If $k=1$, $\Gamma_k$ is contractible in itself and $p>{2n\over
  n-2}$, we prove that there exists $\bar\e>0$ such that, for all
$\e\in(0,\bar\e)$, the Dirichlet problem (\ref{p}) with
$\Omega=T_\e(\Gamma_k)$ does not have any nontrivial solution  (this
nonexistence result follows, as a particular case, from Theorem
\ref{T2.1}).

\no Let us point out that, if $k=1$ but $\Gamma_k$ is noncontractible
in itself or if $k>1$, a nonexistence result analogous to Theorem
\ref{T2.1} cannot hold under the assumption $p>{2n\over  n-2}$.
In fact, the method we use in Theorem \ref{T2.1} fails when $k=1$ and
$\Gamma_k$ is noncontractible because the multipliers to be used in
the Pohozaev type integral identity are not well defined.
Using other multipliers, we obtain a weaker nonexistence result which
holds only when $n\ge 4$ and $p>{2(n-1)\over n-3}$ (it follows from
Theorem \ref{T2.4}).
On the other hand, this weaker result is sharp because, if $\Gamma_k$
is for example a circle of radius $R$ (that is $T_\e(\Gamma_k)$ is a
solid torus), one can easily obtain infinitely many solutions for all
$\e\in(0,R)$ when $n=3$ and $p>1$ or $n\ge 4$ and
$p\in\left(1,{2(n-1)\over n-3}\right)$.

\no Propositions \ref{P3.1}, \ref{P3.2} and \ref{P3.3} give
examples of existence and multiplicity results of positive and sign
changing solutions for some $p\ge{2n\over n-2}$ in tubular domains
$T_\e(\Gamma_k)$ with $k\ge 2$ and $\Gamma_k$ contractible in itself.
This examples explain why Theorem \ref{T2.1} cannot be extended to the
case $k>1$ under the assumption $p>{2n\over n-2}$.

\no However, in the case $k>1$, with $\Gamma_k$ contractible or not,
we prove a weaker nonexistence result (given by Theorem \ref{T3.4})
which holds only when $n>k+2$ and $p>{2(n-k)\over n-k-2}$.

\no Some existence and multiplicity results, when $n\le k+2$ or
$n>k+2$ and $p<{2(n-k)\over n-k-2}$, in tubular domains
$T_\e(\Gamma_k)$ with $k\ge 2$ and $\e$ non necessarily small, show
that also the nonexistence result given by Theorem \ref{T3.4} is
sharp.

\no Finally, let us point out that if in the equation $\Delta
u+f(u)=0$ we replace the Laplace operator $\Delta u$ by the operator
$\div(|D u|^{q-2}Du)$ with $1<q<2$, then critical and supercritical
nonlinearities arise also for $n=2$ and produce analogous nonexistence
results (see \cite{plap,PLap}). 
These results suggest that if $n=2$, $1<q<2$ and $p>{2q\over 2-q}$,
the Pohozaev nonexistence result for  starshaped domains can be
extended to all the contractible domains of $\R^2$ while it is not
possible for example if $n\ge 3$, $q=2$ and $p\ge {2n\over n-2}$
because of Propositions \ref{P3.1}, \ref{P3.2} and \ref{P3.3} (see
Remark \ref{R3.6}).


\sezione{Integral identities and nonexistence results}


In order to obtain nonexistence results for nontrivial solutions of
problem (\ref{*}), we use the Pohozaev type integral identity given in
the following Lemma. 

\begin{lemma}
\label{L2.0}
Let $\Omega$ be a piecewise smooth bounded domain of $\R^n$, $n\ge 3$,
$v=(v_1,\ldots,v_n)\in\cC^1(\overline \Omega,\R^n)$ a vector field in
$\overline\Omega$ and $f$ a continuous function in $\R$. 
Then every solution of problem (\ref{*}) satisfies the integral identity 
\beq
\label{id}
{1\over 2}\int_{\partial\Omega} |Du|^2\, v\cdot\nu\, d\sigma=
\int_\Omega dv[Du]\cdot Du\, dx+
\int_\Omega \div v\, \left(F(u)-{1\over 2}|Du|^2\right)\, dx,
\eeq
where $\nu$ denotes the outward normal to $\partial\Omega$,
$dv[\xi]=\sum\limits_{i=1}^nD_iv\, \xi_i$ $\forall
\xi=(\xi_1,\ldots,\xi_n)\in\R^n$ and $F(t)=\int_0^tf(\tau)d\tau$
$\forall t\in\R$. 
\end{lemma}

\no For the proof it is sufficient to apply the Gauss-Green formula to
the function  $v\cdot Du\, Du$ and argue as in \cite{Po}. 
Notice that the Pohozaev identity is obtained for $v(x)=x$.

\vspace{2mm}

\no Now our aim is to find suitable domains $\Omega$ and vector fields
$v\in\cC^1(\overline \Omega,\R^n)$ 
such that the identity (\ref{id}) can be satisfied only by a trivial
solution of problem (\ref{*}). 

\no In order to construct $\Omega$ and $v$ with this property, let us
consider a curve $\gamma\in\cC^3([a,b],\R^n)$ such that
$\gamma'(t)\neq 0$ $\forall t\in[a,b]$ and
$\gamma(t_1)\neq\gamma(t_2)$ if $t_1\neq t_2$, $t_1,t_2\in [a,b]$. 

\no For all $t\in[a,b]$ and $r>0$, let us set $N(t)=\{\xi\in\R^n\ :\
\xi\cdot\gamma'(t)=0\}$ and $N_r(t)=\{\xi\in N(t)\ :\ |\xi|\le r\}$. 

\no Notice that there exists $\bar\e_1>0$ such that, for all
$\e\in(0,\bar\e_1]$, 
\beq
[\gamma(t_1)+N_\e(t_1)]\cap[\gamma(t_2)+N_\e(t_2)]=\emptyset\quad\mbox{ if
}t_1\neq t_2,\quad  t_1,t_2\in [a,b]. 
\eeq
For all $\e\in(0,\bar\e_1)$ let us consider the open, piecewise
smooth, bounded domain $T^\gamma_\e$ defined by  
\beq
\label{Te}
T^\gamma_\e=\bigcup_{t\in(a,b)} [\gamma(t)+N_\e(t)].
\eeq
Then, the following nonexistence result holds for the nontrivial
solutions in the domain $\Omega=T^\gamma_\e$. 

\begin{teo}
\label{T2.1}
Assume the continuous function $f$ satisfies the condition 
\beq
\label{f} 
tf(t)\ge p\int_0^tf(\tau)d\tau\ge 0\qquad\forall t\in\R
\eeq
for a suitable $p>{2n\over n-2}$.
Then, there exists $\bar \e>0$ such that for all $\e\in(0,\bar\e)$ the
Dirichlet problem (\ref{*}) has only the trivial solution $u\equiv 0$
in the domain $\Omega=T^\gamma_\e$. 
\end{teo}

\no It is clear that condition (\ref{f}) implies $f(0)=0$, so the
function $u\equiv 0$ in $T^\gamma_\e$ is a trivial solution $\forall
\e\in(0,\bar\e_1)$. 

\no In order to prove that it is the unique solution for $\e$ small
enough, we need some preliminary results. 

\no Notice that if $\e\in(0,\bar\e_1)$, the following property holds: for
all $x\in T^\gamma_\e$ there exists a unique $t(x)\in(a,b)$ such that
$\dist(x,\Gamma)=|x-\gamma(t(x))|$, where 
\beq
\Gamma=\{\gamma(t)\ :\ t \in [a,b]\}.
\eeq
If we set $\xi(x)=x-\gamma(t(x))$, we have
$\xi(x)\cdot\gamma'(t(x))=0$ $\forall x\in T^\gamma_\e$. 
Therefore, for all $y\in T^\gamma_\e$ there exists a unique pair $
(t(x),\xi(x))$ such that  $t(x)\in [a,b]$, $\xi(x)\in N_\e(t(x))$ and
$x=\gamma(t(x))+\xi(x)$. 

\no Without any loss of generality, we can assume in addition that
$a\le 0\le b$ and $|\gamma'(t)|=1$ $\forall t\in [a,b]$. 

\no For all $\xi\in N_\e(0)$ let us consider the function $\tau\mapsto
x(\xi,\tau)$ which solves the Cauchy problem 
\beq
\left\{
\begin{array}{l}
\vspace{2mm}
{\partial x\over\partial\tau}=\gamma'(t(x))\\
x(\xi,0)=\gamma(0)+\xi.
\end{array}
\right.
\eeq
Notice that $\dist(x(\xi,\tau),\gamma)=|\xi|$ $\forall \tau\in[a,b]$.
Moreover, for all $\xi\in N_\e(0)$, the function $\tau\mapsto t(x(\xi,\tau))$ is increasing.
As a consequence, we can consider the inverse function $t\mapsto \tau(\xi,t)$ which satisfies $t(x(\xi,\tau(\xi,t)))=t$ $\forall t\in [a,b]$.

\no Notice that $\tau(\xi,0)=0$ $\forall \xi\in N_\e(0)$ because $t(x(\xi,0))=0$.
For all $\xi\in N_\e(0)$, let us set $\psi(\xi,t)=x(\xi,\tau(\xi,t))-\gamma(t)$.
Then, $\psi(\xi,t)\in N_\e(t)$ and $|\psi(\xi,t)|=|\xi|$ $\forall
\xi\in N_\e(0)$ $\forall t\in [a,b]$.
Moreover, for all $x\in T^\gamma_\e$ there exists a unique $\xi\in
N_\e(0)$ such that $\xi(x)=\psi(\xi,t(x))$ and the function
$\xi\mapsto\psi(\xi,t)$ is a one to one function between $N_\e(0)$ and
$N_\e(t)$, satisfying $|\psi(\xi_1,t)-\psi(\xi_2,t)|=|\xi_1-\xi_2|$
$\forall\xi_1,\xi_2\in N_\e(0)$, $\forall t\in[a,b]$.

\no Now, let us consider the vector field $v$ defined by
\beq
\label{v}
v(\gamma(t)+\psi(\xi,t))=t\gamma'(t)[1-\psi(\xi,t)\cdot\gamma''(t)]+\psi(\xi,t)\qquad\forall
t\in(a,b),\ \forall\xi\in N_\e(0).
\eeq
Since $\gamma\in\cC^3([a,b],\R^n)$, we have $v\in\cC^1(\overline
T^\gamma_\e,\R^n)$, so the integral identity (\ref{id}) holds.

\no In the following lemma we estabilish some properties of the vector
field $v$.

\begin{lemma}
\label{L2.2}
In the domain $\overline T^\gamma_{\bar\e_1}$, let us
consider the vector field $v\in \cC^1(\overline T^\gamma_{\bar\e_1},\R^n)$
defined in (\ref{v}).
Then we have
\begin{itemize}
\item[(a)] $v\cdot\nu>0$ on $\partial \overline T^\gamma_\e$
  $\forall\e\in(0,\bar\e_1)$,
\item[(b)] $\lim\limits_{\e\to 0}\sup\{|n-\div v(x)|\ :\ x\in
  T^\gamma_\e\}=0$,
\item[(c)] $\lim\limits_{\e\to 0}\sup\{|1- d v (x)[\eta]\cdot \eta|\ :\ x\in
  T^\gamma_\e,\ \eta\in\R^n,\ |\eta|=1\}=0$.
\end{itemize}
\end{lemma}

\proof Taking into account the choice of $\bar\e_1$, since we are
assuming $|\gamma'(t)|=1$ $\forall t\in[a,b]$, we have
$[1-\psi(\xi,t)\cdot\gamma''(t)]\ge 0$ $\forall t\in [a,b]$.
Therefore, since we are also assuming $a\le 0\le b$, property $(a)$ is
a direct consequence of the definition of $T^\gamma_\e$ and $v$.

\no In order to prove $(b)$, notice that, since $v\in \cC^1(\overline
T^\gamma_\e,\R^n)$ $\forall \e\in(0,\bar\e_1)$, there exist
$t_\e\in[a,b]$ and $\xi_\e\in\overline{N_\e(0)}$ such that
\beq
|n-\div v(\gamma(t_\e)+\psi(\xi_\e,t_\e))|=\max\{|n-\div v(x)|\ :\ x\in
  \overline T^\gamma_\e\}\quad\forall\e\in(0,\bar\e_1).
\eeq
When $\e\to 0$, we obtain (up to a subsequence) $t_\e\to t_0$ for a
suitable $t_0\in[a,b]$ while $\xi_\e\to 0$ (because $|\xi_\e|\le\e$)
and, as a consequence, also $\psi(\xi_\e,t_\e)\to 0$ (because
$|\psi(\xi_\e,t_\e)|=|\xi_\e|$).
Therefore we get 
\beq
\lim_{\e\to 0}\max\{|n-\div v(x)|\ :\ x\in \overline
  T^\gamma_\e\}=|n-\div v(\gamma(t_0))|. 
\eeq
Now, notice that
\beq
dv(\gamma(t_0))[\gamma'(t_0)]=\gamma'(t_0)+t_0\gamma''(t_0) 
\eeq
and 
\beq  
dv(\gamma(t_0))[\psi]=-t_0[\psi\cdot\gamma''(t_0)]\gamma'(t_0)+\psi\qquad\forall\psi\in
N(t_0).
\eeq
It follows that $\div v(\gamma(t_0))=n$, so property $(b)$ holds.

\no In a similar way we can prove property $(c)$.
In fact, since $v\in\cC^1(\overline
T^\gamma_\e,\R^n)$ $\forall \e\in(0,\bar\e_1)$, there exist
$\bar t_\e\in[a,b]$, $\bar \xi_\e\in\overline{N_\e(0)}$ and $\bar\eta
_\e\in \R^n$ such that $|\bar\eta_\e|=1$ and
\beq
|1-dv(\gamma(\bar t_\e)+\psi(\bar\xi_\e,\bar
t_\e))[\bar\eta_\e]\cdot\bar\eta_\e|=
\max 
\{|1- d v (x)[\eta]\cdot \eta|\ :\ x\in
  \overline T^\gamma_\e,\ \eta\in\R^n,\ |\eta|=1\}.
\eeq
Since $|\psi(\bar\xi_\e,\bar t_\e)|=|\bar\xi_\e|\le\e$
$\forall\e\in(0,\bar\e_1)$, we have $\lim\limits_{\e\to
  0}\psi(\bar\xi_\e,\bar t_\e)=0$.
Moreover, there exist $\bar t_0\in  [a,b]$ and $\bar\eta_0\in \R^n$
such that (up to a subsequence) $\bar t_\e\to \bar t_0$ and
$\bar\eta_\e\to\bar\eta_0$ as $\e\to 0$.
It follows that
\beq
\lim_{\e\to 0}\max 
\{|1- d v (x)[\eta]\cdot \eta|\ :\ x\in
  \overline T^\gamma_\e,\ \eta\in\R^n,\ |\eta|=1\}=|1- d v
  (\gamma(\bar t_0))[\bar\eta_0]\cdot \bar\eta_0|.
\eeq
Now, let us set $\bar\psi_0=\bar\eta_0-\bar\eta_0\cdot\gamma'(\bar
t_0)\, \gamma'(\bar t_0)$ and notice that $\bar\psi_0\in N(\bar t_0)$.
Therefore we have 
\beq
dv(\gamma(\bar t_0))[\bar\psi_0]=\bar\psi_0-\bar t_0\,
\bar\psi_0\cdot\gamma ''(\bar t_0)\, \gamma'(\bar t_0).
\eeq
Thus, since
\beq
dv(\gamma(\bar t_0))[\gamma'(\bar t_0)]=\gamma'(\bar
t_0)+\bar t_0\gamma''(\bar t_0)
\eeq
and $\gamma'(\bar t_0)\cdot\gamma''(\bar t_0)=0$, we obtain
\begin{eqnarray}
\nonumber
dv(\gamma(\bar t_0))[\bar\eta_0]\cdot\bar\eta_0 
& = &
dv(\gamma(\bar
t_0))[\bar\eta_0\cdot\gamma'(\bar t_0)\, \gamma'(\bar t_0)+\bar
\psi_0]\cdot
(\bar\eta_0\cdot\gamma'(\bar t_0)\, \gamma'(\bar t_0)+\bar
\psi_0)
\\ \nonumber
& = &
\left\{ \bar\eta_0\cdot\gamma'(\bar t_0)\, [\gamma'(\bar t_0)+\bar
  t_0\gamma''(\bar t_0)]+\bar\psi_0-\bar t_0\,
  \bar\psi_0\cdot\gamma''(\bar t_0)\, \gamma'(\bar t_0)\right\}
\\
& & \cdot (\bar\eta_0\cdot\gamma'(\bar t_0)\, \gamma'(\bar t_0)+\bar\psi_0)
\\ \nonumber
&= &[\bar\eta_0\cdot\gamma'(\bar
t_0)]^2+|\bar\psi_0|^2=|\bar\eta_0|^2=1,
\end{eqnarray}
which implies property $(c)$.

\qed

\begin{cor}
\label{C2.3}
Let $f$ and $F$ be as in Lemma \ref{L2.0}.
Let $T^\gamma_\e$ and $v\in \cC^1(\overline T^\gamma_\e,\R^n)$ be as
in Lemma \ref{L2.2}.
Then, every solution $u_\e$ of the Dirichlet problem (\ref{*}) in
$\Omega =T^\gamma_\e$ satisfies the inequality
\beq
0\le \left[ 1-{n\over
    2}+\mu(\e)\right]\int_{T^\gamma_\e}|Du_\e|^2dx+\int_{T^\gamma_\e}(\div
v)F(u_\e)dx,
\eeq
where $\mu(\e)\to 0$ as $\e\to 0$.
\end{cor}

\no The proof follows directly from Lemmas \ref{L2.0} and \ref{L2.2}.

\vspace{2mm}

{\mbox {{\underline {\sf Proof of Theorem \ref{T2.1}}} \hspace{2mm}}}
In order to prove that the trivial solution $u\equiv 0$ in
$T^\gamma_\e$ is the unique solution for $\e$ small enough, for every
$\e\in(0,\bar\e_1]$, let us consider a solution $u_\e$ of problem (\ref{*}) in
$\Omega=T^\gamma_\e$.  
Taking into account Lemma \ref{L2.0} and condition (\ref{f}), from
Lemma \ref{L2.2} and Corollary \ref{C2.3} we obtain 
\beq
0\le\left[1-{n\over
    2}+\mu(\e)\right]\int_{T^\gamma_\e}|Du_\e|^2dx+[n+\bar\mu(\e)]{1\over
  p}\int_{T^\gamma_\e}u_\e f(u_\e)dx,
\eeq
where $\bar\mu(\e)\to 0$ as $\e\to 0$.
On the other hand, since $u_\e$ is a solution of problem (\ref{*}) in
$\Omega= T^\gamma_\e$, we have
\beq
\int_{T^\gamma_\e}u_\e f(u_\e)\, dx=\int_{T^\gamma_\e}|Du_\e|^2dx.
\eeq
Therefore we obtain
\beq
0\le\left[1-{n\over 2}+{n\over
    p}+\mu(\e)+\bar\mu(\e)\right]\int_{T^\gamma_\e}|Du_\e|^2dx.
\eeq
Since $1-{n\over 2}+{n\over p}<0$ for $p>{2n\over n-2}$, there exists
$\bar \e\in (0,\bar\e_1)$ such that $1-{n\over 2}+{n\over
  p}+\mu(\e)+\bar\mu(\e)<0$ $\forall \e\in(0,\bar\e)$.
Therefore, for all $\e\in(0,\bar\e)$, we must have
$\int_{T^\gamma_\e}|Du_\e|^2dx=0$ which implies $u_\e\equiv 0$ in
$T^\gamma_\e$ and completes the proof.

\qed

\no Notice that if,  instead of the vector field $v$ defined in (\ref{v}), we
consider the vector field $\tilde v$ defined by  
\beq 
\tilde v(\gamma(t)+\psi(\xi,t))=\psi(\xi,t)\qquad\forall t\in
(a,b),\quad\forall \xi\in \overline{N(0)},
\eeq
we obtain a nonexistence  result for $n\ge 4$ and $p>{2(n-1)\over
  n-3}$ (the critical Sobolev exponent in dimension $n-1$, which is
greater than $2n\over n-2$). 

\no Let us point out that the vector field $\tilde v$ is well defined also
when $\gamma$ is a smooth circuit, that is $\gamma(a)=\gamma(b)$ and
$\Omega$ is the interior of $\overline T^\gamma_\e$. 
Therefore, also in these domains we can prove nonexistence results for
$n\ge 4$ and $p>{2(n-1)\over n-3}$, see Theorem \ref{T2.4}. 
On the contrary, in these domains the vector field $v$ could not be
well defined because  
\beq
v(\gamma(a)+\psi(\xi,a))\neq v(\gamma(b)+\psi(\xi,b))\qquad\forall
\xi\in N_\e(0), 
\eeq
while $\gamma(a)+\psi(\xi,a)=\gamma(b)+\psi(\xi,b)$ when
$\gamma(a)=\gamma(b)$ and $\gamma'(a)=\gamma'(b)$. 

\no On the other hand, in these domains one cannot expect to obtain
nonexistence results for $p>{2n\over n-2}$ since it is possible that
there exist nontrivial solutions when $n\ge 4$ and ${2n\over
  n-2}<p<{2(n-1)\over n-3}$ while they do not exist for
$p\ge{2(n-1)\over n-3}$, which happens for example in the case of a solid
torus (see \cite{Pjfa93,Pdie95,MPcras02}). 

\no In next theorem we consider the case where $\Omega$ is a tubular
domain near a circuit, $n\ge 4$ and condition (\ref{f}) holds with
$p>{2(n-1)\over n-3}$ (see Theorem \ref{T3.4} for an extension to more
general tubular domains). 

\begin{teo}
\label{T2.4}
Assume that $\tilde\gamma:[a,b]\to\R^n$ is a smooth curve which
satisfies $\tilde\gamma'(t)\neq 0$ $\forall t\in[a,b]$,
$\tilde\gamma(a)=\tilde\gamma(b)$,
$\tilde\gamma'(a)=\tilde\gamma'(b)$,
$\tilde\gamma(t_1)\neq\tilde\gamma(t_2)$ if $t_1,t_2\in(a,b)$ and
$t_1\neq t_2$. 
Let us set
\beq
\widetilde\Gamma=\{\tilde\gamma(t)\ :\ t\in[a,b]\}
\quad\mbox{ and }\quad
\widetilde T_\e(\widetilde\Gamma)=\{x\in\R^n\ :\
\dist(x,\widetilde\Gamma)<\e\}\quad\forall\e>0. 
\eeq
Moreover assume that $n\ge 4$ and condition (\ref{f}) holds with $p>{2(n-1)\over n-3}$.

\no Then there exists $\tilde\e>0$ such that, for all
$\e\in(0,\tilde\e)$, the Dirichlet problem (\ref{*}) has only the
trivial solution $u\equiv 0$ in the smooth bounded domain
$\Omega=\widetilde T_\e(\widetilde\Gamma)$. 
\end{teo}

\proof
First notice that there exists $\bar\e_1>0$ such that for all
$\e\in(0,\bar\e_1)$ and $x\in\widetilde  T_\e(\widetilde\Gamma)$ there
exists a unique $y\in\tilde\gamma$ such that $\dist(x,\widetilde
\Gamma)=|x-y|$. 
Let us denote this $y$ by $p(x)$ and consider in $\widetilde
T_\e(\widetilde\Gamma)$ the vector field $\tilde v$ defined by $\tilde
v(x)=x-p(x)$. 

\no One can verify by direct computation that 
\beq
d\tilde v(\tilde\gamma(t))[\tilde\gamma'(t)]=0,\quad d\tilde
v(\tilde\gamma(t))[\psi]=\psi\qquad\forall t\in[a,b],\
\forall\psi\in\R^n\ \mbox{ such that }\ \psi\cdot\tilde\gamma'(t)=0 
\eeq
and, as a consequence, 
\beq
\label{e1}
\div\tilde v(\tilde\gamma(t))=n-1\qquad\forall t\in[a,b]
\eeq
\beq
\label{e2}
d\tilde
v(\tilde\gamma(t))[\eta]\cdot\eta=|\eta|^2-{[\eta\cdot\tilde\gamma'(t)]^2\over
  |\tilde\gamma'(t)|^2}\qquad \forall t\in[a,b],\ \forall
\eta\in\R^n. 
\eeq
It follows that 
\beq
\label{e1'}
\lim_{\e\to 0}\sup\{|n-1-\div\tilde v(x)|\ :\ x\in\widetilde  T_\e(\widetilde\Gamma)\}=0
\eeq
as one can easily obtain from (\ref{e1}) arguing as in the proof of
assertion $(b)$ of Lemma \ref{L2.2}. 
Moreover, from (\ref{e2}) we obtain
\beq
\label{e2'}
\lim_{\e\to 0}\sup\{d\tilde v(x)[\eta]\cdot\eta\ :\ x\in\widetilde
T_\e(\widetilde\Gamma),\ \eta\in \R^n,\ |\eta|=1\}=1. 
\eeq 
In fact, for all $\e\in(0,\bar\e_1)$, choose $x_\e\in \widetilde
T_\e(\widetilde\Gamma)$ and $\eta_\e\in\R^n$ such that $|\eta_\e|=1$
and $s_\e-\e\le d\tilde v(x_\e)[\eta_\e]\cdot\eta_\e$ where 
\beq
s_\e=\sup\{d\tilde v(x)[\eta]\cdot\eta\ :\ x\in\widetilde
T_\e(\widetilde\Gamma),\ \eta\in \R^n,\ |\eta|=1\}. 
\eeq
Since $\dist(x_\e,\widetilde\Gamma)\to 0$ as $\e\to 0$, and
$\widetilde\Gamma$ is a compact manifold, from (\ref{e2}) we infer
that $\limsup\limits_{\e\to 0} s_\e\le 1$. 
On the other hand, (\ref{e2}) implies $s_\e\ge 1$
$\forall\e\in(0,\bar\e_1)$, so (\ref{e2'}) is proved. 

\no Furthermore, one can easily verify that $\tilde v\cdot\nu>0$ on
$\partial \widetilde  T_\e(\widetilde\Gamma)$
$\forall\e\in(0,\bar\e_1)$. 
Thus, taking also into account condition (\ref{f}), from Lemma
\ref{L2.0} we infer that every solution $\tilde u_\e$ of problem
(\ref{*}) in the domain $\widetilde  T_\e(\widetilde\Gamma)$ satisfies 
\beq
0\le \left[ 1-{n-1\over 2}+\tilde\mu(\e)\right]\int_{\widetilde
  T_\e(\widetilde\Gamma)}|D\tilde u_\e|^2dx+\left[{n-1\over
    p}+\tilde\mu(\e)\right]\int_{\widetilde
  T_\e(\widetilde\Gamma)}\tilde u_\e f(\tilde u_\e)\, dx, 
\eeq
where $\tilde\mu(\e)\to 0$ as $\e\to 0$.
Since
\beq
\int_{\widetilde  T_\e(\widetilde\Gamma)}\tilde u_\e f(\tilde
u_\e)\, dx=\int_{\widetilde  T_\e(\widetilde\Gamma)}|D\tilde u_\e|^2dx 
\eeq
(because $\tilde u_\e$ solves problem (\ref{*}) in $\widetilde
T_\e(\widetilde\Gamma)$) we obtain 
\beq
\label{le}
0\le\left[1-{n-1\over 2}+{n-1\over p}+2\tilde
  \mu(\e)\right]\int_{\widetilde  T_\e(\widetilde\Gamma)}|D\tilde
u_\e|^2dx 
\eeq
where $1-{n-1\over 2}+{n-1\over p}<0$ because $n\ge 4$ and $p>{2(n-1)\over n-3}$.
Therefore, there exists $\tilde \e\in(0,\bar\e_1)$ such that for all
$\e\in(0,\tilde\e)$ (\ref{le}) implies $\tilde u_\e\equiv 0$ in  
$\widetilde  T_\e(\widetilde\Gamma)$.
So the proof is complete.

\qed


\sezione{Tubular domains of higher dimension and final remarks}


The nonexistence results presented in Section 2 are concerned with
domains $\Omega$ which are thin neighbourhoods of 1-dimensional
manifolds (with boundary and contractible in Theorem \ref{T2.1},
without boundary and noncontractible in Theorem  \ref{T2.4}).  
In this section we consider the case where $\Omega$ is a thin
neighbourhood of  $k$-dimensional smooth, compact manifold $\Gamma_k$
with $k>1$. 

\no If $\Gamma_k$ is a submanifold of $\R^n$ with $n>k$, for all
$x\in\Gamma_k$ we set $N(x)=T^\perp(x)$  and $N_\e(x)=\{x\in N(x)\ :\
|x|<\e\}$, where $T(x)$ is the tangent space to $\Gamma_k$ in $x$ and
$N(x)$ is the normal space. 
Since $\Gamma_k$ is a compact smooth submanifold, there exists
$\bar\e_1>0$ such that, for all $\e\in(0,\bar\e_1]$,  we have
$[x_1+N_\e(x_1)]\cap [x_2+N_\e(x_2)]=\emptyset$ for all $x_1$ and $x_2$
in $\Gamma_k$ such that $x_1\neq x_2$. 
Then, for all $\e\in(0,\bar\e_1)$, we consider the piecewise smooth,
bounded domain $T_\e(\Gamma_k)$ defined as the interior of the set 
$\cup_{x\in\Gamma_k}[x+N_\e(x)]$ (we say that $T_\e(\Gamma_k)$ is the
tubular domain with thickness $\e$ and center $\Gamma_k$). 
Our aim is to study existence and nonexistence of nontrivial solutions
of problem (\ref{*}) in the domain $\Omega=T_\e(\Gamma_k)$. 

\no Let us point out that when $k>1$ one cannot prove a theorem
analogous to Theorem \ref{T2.1}. 
In fact, if $\Gamma_k$ is a $k$-dimensional manifold contractible in
itself and $k>1$, one cannot obtain nonexistence results for
nontrivial solutions of problem (\ref{*}) in the domain
$\Omega=T_\e(\Gamma_k)$ under the assumption that condition (\ref{f})
holds with $p>{2n\over n-2}$ as in Theorem \ref{T2.1}. 
The reason is explained by the following examples where existence results hold.
\begin{ex}
\label{E3.0}
{\em
For all $n\ge k+1$, let us consider the function
$\gamma_k:\R^k\to\R^n$ defined as follows: 
\beq
\begin{array}{rcll}
\gamma_{k,i}(x_1,\ldots,x_k)&=& {2x_i\over |x|^2+1}&\qquad \mbox{ for }i=1,\ldots,k
\\
\gamma_{k,k+1}(x_1,\ldots,x_k)&=& {|x|^2-1\over |x|^2+1} &
\\
\gamma_{k,i}(x_1,\ldots,x_k)&=& 0 &\qquad\mbox{ for }i=k+2,\ldots,n
\end{array}
\eeq
($\gamma_k$ is the stereographic projection of $\R^k$ on a
$k$-dimensional sphere of $\R^n$). 

\no Moreover, for all $r>0$, let us set $\Gamma^r_k=\{\gamma_k(x)\ :\ x\in\R^k,\ |x|<r\}$.
}
\end{ex}

\no Then one can easily verify that the domain $T_\e(\Gamma^r_k)$ is
contractible in itself for all $r>0$ and $\e\in(0,1)$. 
Moreover, the following propositions hold.

\begin{prop}
\label{P3.1}
Let $k\ge 2$ and $n\ge k+1$.
Assume that $f(t)=|t|^{p-2}t$ with $p\ge{2n\over n-2}$ and that
  $p<{2(n-k+1)\over n-k-1}$ if $n>k+1$. 

\no Then, there exists $\bar r>0$ such that if $r>\bar r$ and
$\e\in(0,1)$, problem (\ref{*}) in the domain
$\Omega=T_\e(\Gamma^r_k)$ has positive and sign changing solutions;
moreover, under the additional assumption $p>{2n\over n-2}$, for all
$\e\in(0,1)$ the number of solutions tend to infinity as $r\to\infty$. 
\end{prop}

\no For the proof it suffices to look for solutions having radial
symmetry with respect to the first $k$-variables and argue as in
\cite{Pmm89,Pd98,P4092,Ptmna96,P94,MPaihp06,MPcvpde06,MPaihp04}. 

\begin{prop}
\label{P3.2}
Let $k\ge 2$, $n\ge k+1$, $r>1$, $\e\in(0,1)$.
Moreover, assume that $f(t)=|t|^{p-2}t$ $\forall t\in\R$.
Then, there exists $\bar p>{2n\over n-2}$ such that, if $n=k+1$ and
$p\ge \bar p$ or if $n>k+1$ and $p\in\left[\bar p,{2(n-k+1)\over
    n-k-1}\right)$, problem (\ref{*}) with $\Omega=T_\e(\gamma^r_k)$
has solution. 
\end{prop}

\no The proof can be carried out arguing for example as in
\cite{MPcvpde06} in order to obtain solutions having radial symmetry
with respect to the first $k$ variables. 

\begin{prop}
\label{P3.3}
Let $k\ge 2$, $n\ge k+1$, $r>1$, $\e\in(0,1)$ and assume that
$f(t)=|t|^{p-2}t$ $\forall t\in\R$. 
Then, there exists $\tilde p>{2n\over n-2}$ such that problem
(\ref{*}) with $\Omega=T_\e(\Gamma^r_k)$ has positive solutions for
all $p\in\left({2n\over n-2} ,\tilde p\right)$. 
Moreover, the number of solutions tends to infinity as $p\to{2n\over
  n-2}$. 
\end{prop}

\no The proof is based on a Lyapunov-Schmidt type finite dimensional
reduction method as in \cite{MPcras02},\cite{MPaihp04}, etc.  

\no Thus, while Theorem \ref{T2.1} gives a nonexistence result for all
$p>{2n\over n-2}$ when $k=1$, $\Gamma_k$ is contractible in itself and
$\Omega$ is a thin tubular domain centered in $\Gamma_k$, Propositions
\ref{P3.1}, \ref{P3.2} and \ref{P3.3} give examples of existence
results for some $p>{2n\over n-2}$ when $\Omega$ is a tubular domain
centered in a suitable $k$-dimensional manifold $\Gamma^r_k$,
contractible in itself but with $k\ge 2$. 
In this sense we mean that Theorem \ref{T2.1} cannot be extended to
the case $k\ge 2$ (see also Remark \ref{R3.5} for more details about
the differences between the cases $k=1$ and $k>1$). 
 
\no However, notice that a weaker nonexistence result holds for all
$k\ge 1$ (even if $\Gamma_k$ is noncontractible in itself) when
$n>k+2$ and $p>{2(n-k)\over n-k-2}$, as we prove in the following
Theorem \ref{T3.4}. 

\no If $n\le k+2$ or $n>k+2$ and $p<{2(n-k)\over n-k-2}$, the
existence of nontrivial solutions can be proved even if $\Omega$ is a
tubular domain $T_\e(\Gamma_k)$ with $\e$ not necessarily small: for
example, if $\Gamma_k$ is a $k$-dimensional sphere, we can look for
solutions with radial symmetry with respect to $k+1$ variables, so we
obtain infinitely many solutions for all $\e\in(0,R)$ where $R$ is the
radius of the sphere. 

\begin{teo}
\label{T3.4}
Let $k\ge 1$, $n>k+2$ and assume that $\Gamma_k$ is a $k$-dimensional,
compact, smooth submanifold of $\R^n$.
Moreover, assume that condition (\ref{f}) holds with $p>{2(n-k)\over
  n-k-2}$.

\no Then, there exists $\bar\e>0$ such that, for all
$\e\in(0,\bar\e)$, the Dirichlet problem (\ref{*}) has only the
trivial solution $u\equiv 0$ on the tubular domain
$\Omega=T_\e(\Gamma_k)$.
\end{teo}

\proof 
Taking into account the definition of the tubular domain
$T_\e(\Gamma_k)$, for all $\e\in(0,\bar\e_1)$ and $x\in T_\e(\Gamma_k)$
there exists a unique $y\in\Gamma_k$ such that $x\in y+N_\e(y)$.
Then, denote this $y$ by $p_k(x)$ and set $v_k(x)=x-p_k(x)$ $\forall
x\in T_\e(\Gamma_k)$.
One can easily verify that the vector field $v_k$ satisfies
$v_k\cdot \nu\ge 0$ on $\partial T_\e(\Gamma_k)$
$\forall\e\in(0,\bar\e_1)$.

\no Therefore, from Lemma \ref{L2.0} we infer that every solution
$u_\e$ of problem (\ref{*}) in $T_\e(\Gamma_k)$ satisfies
\beq
\label{lek}
0\le\int_{T_\e(\Gamma_k)}dv_k[Du_\e]\cdot Du_\e\,
dx+\int_{T_\e(\Gamma_k)}\div v_k\left(F(u_\e)-{1\over 2}|D
  u_\e|^2\right)dx.
\eeq
Notice that 
\beq
dv_k(x)[\phi]=0,\qquad dv_k(x)[\psi]=\psi\qquad \forall x\in\Gamma_k,\
\forall\phi\in T(x),\ \forall\psi\in N(x)
\eeq
as one can verify by direct computation.

\no As a consequence we obtain
\beq
\div v_k(x)=n-k,\qquad dv_k(x)[\phi+\psi]\cdot (\phi+\psi)=|\psi|^2
\qquad \forall x\in\Gamma_k,\
\forall\phi\in T(x),\ \forall\psi\in N(x).
\eeq
Since $\Gamma_k$ is a compact manifold, it follows that 
\beq
\lim_{\e\to 0}\sup\{|n-k-\div v_k(x)|\ :\ x\in T_\e(\Gamma_k)\}=0
\eeq
and
\beq
\lim_{\e\to 0}\sup\{dv_k(x)[\eta]\cdot\eta\ :\ x\in T_\e(\Gamma_k),\
\eta\in\R^n,\ |\eta|=1\}=1
\eeq
as one can infer arguing as in the proof of Theorem \ref{T2.4}.

\no Thus, taking also into account that
\beq
\int_{T_\e(\Gamma_k)} u_\e f(u_\e)\, dx=
\int_{T_\e(\Gamma_k)}|Du_\e|^2dx,
\eeq
from condition (\ref{f}) we infer that
\beq
0\le \left[1-{n-k\over 2}+{n-k\over
    p}+\mu_k(\e)\right]\int_{T_\e(\Gamma_k)} |Du_\e|^2dx
\eeq
where $\mu_k(\e)\to 0$ as $\e\to 0$.
Since $1-{n-k\over 2}+{n-k\over p}<0$ (because $n>k+2$ and
$p>{2(n-k)\over n-k-2}$), it follows that there exists
  $\bar\e\in(0,\bar\e_1)$ such that, for all $\e\in(0,\bar\e)$ we have
  $u_\e\equiv 0$ in $T_\e(\Gamma_k)$, so the problem has only the
  trivial solution $u\equiv 0$.

\qed

\begin{rem}
\label{R3.5}
{\em
Proposition \ref{P3.1}, as well as the results reported in
\cite{Pmm89,Pd98,P4092,Ptmna96,P94,MPaihp06,MPcvpde06,MPaihp04},
suggest that the existence of nontrivial solutions is related to the
property that the domain $\Omega$ is obtained by removing a subset of
small capacity from a domain having a different $k$-dimensional
homology group with $k\ge 2$.

\no For example, in the case of domains with small holes, every hole
has small capacity and changes the $(n-1)$-dimensional homology group.

\no In the case of tubular domains $T_\e(\Gamma^r_k)$, the existence
results for $k\ge 2$ and $r$ large enough given by Proposition
\ref{P3.1} is related to the fact that $\Gamma^r_k$ tends to a
$k$-dimensional sphere $S_k$ as $r\to\infty$, the capacity of 
$T_\e(S_k)\setminus T_\e(\Gamma^r_k)$ tends to 0 as $r\to\infty$ and the
domains $T_\e(S_k)$ and  $T_\e(\Gamma^r_k)$ have different
$k$-dimensional homology group. 

\no On the contrary, when $k=1$, the capacity of $T_\e(S_1)\setminus
T_\e(\Gamma^r_1)$ does not tend to 0 as $r\to\infty$.
This fact explains the nonexistence result given by Theorem \ref{T2.1}
in the case of the domains $T_\e(\Gamma^r_1)$, when $\e$ is small
enough, for all $r>0$.
}\end{rem}

\begin{rem}
\label{R3.6}
{\em
If $n=2$ we do not have critical or supercritical phenomena for the
Laplace operator.
But, if we replace it by the $q$-Laplace operator, this phenomena
arise and may produce nonexistence results for nontrivial solutions.
For example, if we consider the Dirichlet problem
\beq
\label{q}
\div(|Du|^{q-2}Du)+|u|^{p-2}u=0\quad\mbox{ in }\Omega,\qquad
u=0\quad\mbox{ on }\partial \Omega
\eeq
where $\Omega$ is a bounded domain of $\R^2$, $1<q<2$, $p\ge {2q\over
  2-q}$, then one can prove nonexistence results in some bounded
contractible domains which can be non starshaped and even arbitrarily
close to noncontractible domains (see \cite{plap,PLap}).
For example, if $\Omega=T_\e(\Gamma^r_1)$, there exists $\bar\e>0$
such that problem (\ref{q}) has only the trivial solution $u\equiv 0$
for all $r>0$ and $\e\in(0,\bar\e)$.

\no The results obtained in \cite{plap,PLap} suggest that the
nonexistence of nontrivial solutions for Dirichlet problem (\ref{q}) 
might be proved in all the contractible domains of $\R^2$ (while it is not
possible for problem (\ref{p}) when $n\ge 3$ and $p\ge {2n\over n-2}$
because of Proposition \ref{P3.1}).
}\end{rem}


{\small {\bf Acknowledgement}. The authors have been supported by the ``Gruppo
Nazionale per l'Analisi Matematica, la Probabilit\`a e le loro
Applicazioni (GNAMPA)'' of the {\em Istituto Nazio\-nale di Alta Matematica
(INdAM)} - Project: Equazioni di Schrodinger nonlineari: soluzioni con
indice di Morse alto o infinito. 

The second author acknowledges also the MIUR Excellence Department
Project awarded to the Department of Mathematics, University of Rome
Tor Vergata, CUP E83C18000100006 
}



\end{document}